\documentclass[11pt]{article}
\usepackage{graphicx}
\usepackage{amssymb}
\usepackage{amsmath,amsthm,amscd,amssymb}
\usepackage{latexsym}
\numberwithin{equation}{section}

\begin{document}
\theoremstyle{plain}
\newtheorem{theorem}{Theorem}[section]
\newtheorem{lemma}[theorem]{Lemma}
\newtheorem{corollary}[theorem]{Corollary}
\newtheorem{proposition}[theorem]{Proposition}
\newtheorem{conjecture}[theorem]{Conjecture}
\newtheorem{criterion}[theorem]{Criterion}
\newtheorem{algorithm}[theorem]{Algorithm}

\theoremstyle{definition}
\newtheorem{definition}[theorem]{Definition}
\newtheorem{condition}[theorem]{Condition}
\newtheorem{problem}[theorem]{Problem}
\newtheorem{example}[theorem]{Example}
\newtheorem{exercise}{Exercise}[section]
\newtheorem{obs}{Observation}

\theoremstyle{remark}
\newtheorem{remark}[theorem]{Remark}
\newtheorem{question}[theorem]{Question}
\newtheorem{note}[theorem]{Note}
\newtheorem{notation}[theorem]{Notation}
\newtheorem{claim}[theorem]{Claim}
\newtheorem{summary}[theorem]{Summary}
\newtheorem{acknowledgment}[theorem]{Acknowledgment}
\newtheorem{case[theorem]}{Case}
\newtheorem{conclusion}[theorem]{Conclusion}

\def \R{{\mathbb R}}
\def \F{{\mathbb F}}
\def \Z{{\mathbb Z}}

\title{{\Large On distance measures for well-distributed sets \footnotetext{The work was partly supported by the grant DMS02-45369 from the National Science
Foundation, the National Science Foundation Focused Research Grant DMS04-56306, and the EPSRC grant
GR/S13682/01.} }\footnotetext{AMS subject classification 42B, 52C, 97C} }

\author{A. Iosevich\thanks{University of Missouri, Columbia MO, 65211
USA, {\tt iosevich@math.missouri.edu}} \ and M. Rudnev\thanks{University of Bristol, Bristol BS8 1TW UK, {\tt
m.rudnev@bris.ac.uk}}}

\maketitle

\begin{abstract} In this paper we investigate the Erd\"os/Falconer distance conjecture for a natural class of sets
statistically, though not necessarily arithmetically, similar to a lattice. We prove a good upper bound for
spherical means that have been classically used to study this problem. We conjecture that a majorant for the
spherical means suffices to prove the distance conjecture(s) in this setting. For a class of non-Euclidean
distances, we show that this generally cannot be achieved, at least in dimension two, by considering integer
point distributions on convex curves and surfaces. In higher dimensions, we link this problem to the question
about the existence of smooth well-curved hypersurfaces that support many integer points.
\end{abstract}


\section{Introduction}

In this paper we study the Erd\"os/Falconer distance problems, introduced in \cite{Er} and \cite{Fa},
respectively, which ask, in discrete and continuous setting, respectively, whether an appropriate size condition
on a subset of the Euclidean space guarantees that the set of pair-wise distances determined by the set is also
suitably large. In both the continuous and discrete analogs of the problem, the integer lattice and its suitable
thickenings and scalings provide the most compelling examples indicating the sharpness of the conjectured
results. This led several authors in recent years to study well-distributed sets, which are statistically,
though not necessarily arithmetically, analogous to the integer lattice. In particular, it is observed in
\cite{IL} and \cite{HI} that in this context the estimate for the continuous problem can be readily converted
into a corresponding estimate for its discrete predecessor. We shall see below that the structure of these sets
lends itself to elegant and relatively straight forward analysis using Fourier analytic methods. We shall also
use arithmetic considerations to indicate limitations of these methods for Euclidean and non-Euclidean metrics,
especially the latter. The paper is concluded with an explicit synthesis of analytic, combinatorial and number
theoretic considerations in the context of incidence problems.

In the continuous setting, the main tool in the study of the distance set problem, pioneered by Mattila
(\cite{Ma}) is the $L^2$ spherical average of the Fourier transform. The two main results of this paper deal
with this fundamental object. In the positive direction we obtain a good upper bound for these averages in the
context of measures obtained by thickening well-distributed sets. These estimates are stronger than the
corresponding results previously obtained by Wolff (\cite{WA}), Erdo\~gan (\cite{Erg}), and others in the
context of general measures.

Our estimates are based on the estimate of a certain exponential sum that naturally appears in our calculations.
Coarse bounds for this sum enable us to match the generally optimal spherical average upper bound of Wolff
(\cite{WI}) in dimension two where we also gain the endpoint. In higher dimensions, they yield better estimates
for the spherical average than the best known general bounds due to Erdo\~gan (\cite{Er}).  As we note above,
our gain, or rather its technical transparency, is due to special features of the well-distributed set-up,
whereas the bounds of Wolff and Erdo\~gan apply to general Frostman measures. We do expect however that our
improved bound in higher dimension should hold in the general case.

The fact that Wolff's general bound in two dimensions cannot be improved is supported by a counterexample of Sj\"olin (\cite{Sj}) although the latter is highly not well-distributed. There is no evidence that in the well-distributed setting, the spherical average does not satisfy sufficient good bounds to imply the Erd\"os distance conjecture, and  we conjecture that this is indeed the case. Our conjecture is supported by the integer lattice case, when the coarse upper bound for the aforementioned sum can be easily refined by using the Poisson summation formula and elementary number theory. We conclude that the well-distributed Erd\"os conjecture may well follow from obtaining sufficiently sharp estimates for the sum in question.

Estimates for the spherical average provide lower bounds for the number of distinct distances regarding the Erd\"os conjecture. To this end, our coarse bound enables us to match earlier results of Moser (\cite{Mo}) Solymosi and Vu (\cite{SV}), and one of the authors (\cite{Io}) obtained by purely combinatorial methods.

Our second result, a lower bound on spherical averages, provides evidence that in a broader setting of non-isotropic distances generated by well-curved smooth convex bodies, majorants for the analogs of the spherical average alone, which are the Fourier $L^2$ averages over dilates of the boundary of the dual body, do not generally imply the Erd\"os conjecture even in the well-distributed case. This makes the Euclidean distance special, as it is in the case of the single distance conjecture in the plane which is generally false for non-Euclidean distances. The counterexample we present is built on the integer lattice in two dimensions and the fact that there exist well-curved domains whose dilated boundary hosts many more integer points then does the circle. In this context, we provide a bound from below on how effective the corresponding majorant for the $L^2$ average over the boundary can be. In the last section of the paper, these results are discussed in some detail in regard to how the distance conjectures can be related to the problem of lattice points distributions on dilates of convex curves and hypersurfaces.

\subsubsection*{Distance conjectures and statement of results}

The integer lattice case prompted Erd\"os to conjecture that for a point set $E \subset {\mathbb R}^d$, its
distance set
\begin{equation}\label{ds} \Delta(E)=\{\|x-y\|: x,y \in E\} \end{equation} is such that
\begin{equation}\label{ec} |\Delta(E)| \geq c(|E|) |E|^{\frac{2}{d}}, \end{equation}
where the ``constant'' $c(|E|)$ does not vanish faster than any negative power of the cardinality
$|E|\rightarrow\infty$.  In the sequel, $\|\cdot\|$ is the Euclidean norm, the notation $|E|$ is used to denote the cardinality of a discrete set or the Lebesgue measure of a continuum set $E\subset \R^d$.

Thickening of the lattice was also used by Falconer to support the conjecture that for a Borel set
$E\subset \R^d$,
\begin{equation}\label{fc}{\rm dim}_H\,E>\frac{d}{2} \Rightarrow |\Delta(E)|>0,
\end{equation}
where ${\rm dim}_H$ is the Hausdorff dimension and $|\cdot|$ is the Lebesgue measure.

\vskip.125in

The Falconer conjecture can be regarded as the ``continuous version'' of the Erd\"os conjecture, though a quantitative link, obtained in \cite{HI} and \cite{IL} is only known in the context of well-distributed sets.

We say that an infinite point set $A \subset {\mathbb R}^d$ is class ${\cal A}$ well-distributed (sometimes also
known as homogeneous, or Delaunay which some authors spell as {\em Delone}) if it is separated in the sense that
for some $c_{\cal A}$, one has $\|a-a'\| \geq c_{\cal A}, \forall \ a,a' \in A: a \neq a',$ as well as any cube
of side length $C_{\cal A}$ has a non-empty intersection with $A$. Constants in the ensuing estimates related to
$A \in {\cal A}$ will depend on $c_{\cal A}$ and $C_{\cal A}$, which are regarded as fixed and denoted by $c$
and $C$, respectively.

For the truncations $A_q= A \cap B(0,q)$ of $A \in {\cal A}$, with $q\gg1$ and $B(x,\delta)$ denoting the Euclidean ball of radius $\delta$ centered at $x$, the Erd\"os conjecture says that
\begin{equation} |\Delta(A_q)| \gtrapprox q^2,\label{wde} \end{equation} where $q$ naturally becomes a parameter, a slow dependence in which may be hidden (as it is in the case $d=2$) in the symbol
$\gtrapprox$.

The notation $X \lesssim Y$ ($X\gtrsim Y$) means that for some large positive $C$ (small positive $c$), one has $X \leq CY$ ($X\geq cY$). Furthermore, $X\approx Y$ if both $X\lesssim Y$ and
$X\gtrsim Y$. The notations $X \lesssim Y$, $X\gtrsim Y$, and $X\approx Y$  also appear in the literature as $X=O(Y)$, $X=\Omega(Y)$, and $X=\Theta(Y)$, respectively, and we use these as
well. Besides, the notations $X \lessapprox Y$ ($X \gtrapprox Y$) mean that the constant hidden in the
$\lesssim$ ($\gtrsim$) symbol may be allowed to grow (decrease) slower than any power of the controlling parameter that is associated with the estimate.

Note that the distance conjecture (\ref{wde}) in the well-distributed setting in $d>2$ would follow from the
case $d=2$ by restricting the set $A_q \in {\mathbb R}^d$ to a ``horizontal" slab of thickness $O(1)$ in
${\mathbb R}^{d-1}$.

To bring Fourier analysis into the problem, given a Delaunay set $A$, one thickens it to create a Cantor-like set $E$ of any dimension $s=d/p,\,p>1$. See, for example, \cite{Fa}, \cite{HI} and \cite{IL}. The first step in the construction is scaling  the truncation $A_q$ into the unit ball and then thickening each point into a ball of radius $\Theta(q^{-p}).$ Let us call the resulting set $E_q$. Then the number of $q^{-p+1}$ separated distances that $A_q$ generates can be estimated by $\Omega(q^{p}|\Delta(E_q)|).$

More precisely, in what follows let $\phi(x), x \in {\mathbb R}^d$  be  a radial test function, whose support is
contained in the unit ball. Suppose that $\phi$ is positive in the interior of its support, $\phi(0)=1$,
$\int\phi=1$, and the Fourier transform $\widehat{\phi}$ is non-negative. Let $\phi_{q^p}(x) = q^{pd} \phi(q^p
x)$ and define
\begin{equation}\label{measure} d\mu_s(x)=\rho_s(x)dx, \mbox{ with } \rho_s(x)=C_A q^{-d} \sum_{a \in q^{-1}A_q} \phi_{q^p}(x-a),\end{equation} where $C_A$ is the normalization constant. Heuristically,
\begin{equation}\label{mu}\rho_s(x) \sim q^{(p-1)d} \sum_{a \in q^{-1}A_q}B_{a,q^{-p}}(x),\end{equation} where the notation $B_{a,q^{-p}}$ for the ball centered at $a$, of radius $q^{-p}$ has been identified with its characteristic function $B_{a,q^{-p}}(x)$.

The Lebesgue measure $|\Delta(E_q)|$ can be bounded from below by methods that have been developed for the
Falconer distance problem and are the main concern of this paper. In general, for a compact Borel $E \subset
{\mathbb R}^d,$ a Borel probability measure $\mu$ supported on $E$ defines automatically the distance measure
$\nu_{\mu}$ as the push-forward of $\mu \times \mu$ under the distance map $E \times E \mapsto
\Delta(E)\subset\R_+$.

Such a set $E$, with  ${\rm dim}_H(E)=\alpha,$ supports for all $s<\alpha$ a Frostman measure $\mu$, so that
\begin{equation} \int_{B(x,\delta)}d\mu \leq C_\mu \delta^s, \forall\,x \mbox{ and all } \delta\ll1,
\label{fm}\end{equation} and
\begin{equation}
I_s(\mu)=\int\int\frac{d\mu_xd\mu_y}{\|x-y\|^s}<\infty.
\label{energy}\end{equation}

If ${\cal M}_s$ is a class of such measures and $\mu\in{\cal M}_s$,
an important sub-problem in the Falconer conjecture is to establish
general asymptotic bounds for the spherical average
\begin{equation}\label{savr} \sigma_\mu(t)\;=\; \int_{S^{d-1}}
{|\widehat{\mu}(t \omega)|}^2 d\omega,
\end{equation}
in the form \begin{equation}
\label{ub}\sigma_\mu(t)\;\leq\;C_{\mu,\beta}
t^{-\beta},\;\;\forall\,t\geq 1.\end{equation}

Given $s$ and $\mu\in{\cal M}_s$, the best known results are as follows: the bound (\ref{ub}) holds for all
\begin{equation}\beta\;\;<\;\;\left\{
\begin{array}{lllll} s,& \mbox{for}& 0<s\leq\frac{d-1}{2}, \\ \hfill \\
\frac{d-1}{2},& \mbox{for}& \frac{d-1}{2}\leq s\leq\frac{d}{2}, \\ \hfill \\
\frac{d+2s-2}{4},& \mbox{for}& \frac{d}{2}\leq s\leq\frac{d+2}{2}, \\ \hfill \\
 s-1,& \mbox{for}& \frac{d+2}{2}\leq s<d.\end{array}\right.\label{gnr}\end{equation}
These results are due to Falconer (\cite{Fa}), Mattila (\cite{Ma}), Sj\"olin (\cite{Sj}), Wolff (\cite{WI},
Erdo\~gan (\cite{Er}), and others, see, for examples, the references contained in \cite{Er}. The crucial
interval of the values of $s$ where one would like to improve over (\ref{gnr}) is for
$s\in[\frac{d}{2},\frac{d+1}{2}]$, and first and foremost at the ``critical'' value $s=\frac{d}{2}.$ We will
develop more background in further chapters of the paper and will now formulate its main results.

\begin{theorem}\label{tone} Let $A\in{\cal A}$ be a well-distributed set, let the measure
$\mu_s$ be defined by (\ref{measure}), with $p \in (1,2]$. Then $\mu_s\in{\cal M}_{s=\frac{d}{p}}$ and for some $\eta(\tau)\leq C_n(1+|\tau|)^{-n},$ for any large $n$, one has
\begin{equation}\label{upper}
\sigma_{\mu_s}(t)\;\lesssim\; q^{-d+1} \eta\left(\frac{t}{q^p}\right).
\end{equation}
In addition, for $s=\frac{d}{2}$ and  $t\gtrsim q$,
\begin{equation}\label{sigmabound}\hspace{-5mm}
\sigma_{\mu_{d/2}}(t)\;\lesssim \;t^{1-d}\Sigma_{d/2}(t),\end{equation} where the quantity $\Sigma_{d/2}(t),$ to be defined explicitly, satisfies the coarse bound
\begin{equation}\label{coarse} \Sigma_{d/2}(t) \lesssim t^{\frac{d-1}{2}}\eta\left(\frac{t}{q^2}\right).
\end{equation}
\end{theorem}
\begin{remark} The bound (\ref{coarse}) for the quantity $\Sigma_{d/2}(t)$ does not appear to be optimal, an one can expect that it can be improved, by the factor of $\sqrt{t}$. This would then imply the Erd\"os conjecture for well-distributed sets. The explicit expression for $\Sigma_{d/2}(t)$ as well as a conditional bound, which, as suggested by the integer lattice example discussed at the final section of the paper, can indeed beat (\ref{sigmabound}) by the factor $\sqrt{t}$ are given by (\ref{sigma}) below.\end{remark}

By (\ref{upper}), with $s\geq\frac{d}{2}$ we have a bound
\begin{equation}\label{cb}\sigma_{\mu_s}(t)\;\lesssim\; t^{-\frac{d-1}{d}s},\;\;\forall\,t\geq 1.\end{equation}
which has the following corollary.
\begin{corollary}\label{corone}
The set $A_q$ determines  $\Omega\left(|A_q|^{\frac{2}{d}-\frac{1}{d^2}}\right)$ distinct distances.
\end{corollary}

It is interesting to broaden the scope of the distance conjectures by generalizing the Euclidean distance
$\|\cdot\|$ as $\|\cdot\|_K$, the Minkowski functional of a strictly convex body $K \subset {\mathbb R}^d,$ with
the smooth boundary $\partial K$. Let ${\cal K}$ be described a class of such bodies, whose volume equals the
volume of the unit ball, and the Gaussian curvature is bounded in some interval $[c_{\cal K},C_{\cal K}]$. The
spherical average generalizes accordingly by replacing the domain of integration $S^{d-1}$ in (\ref{savr}) by
$\partial K$, substituting $S^{d-1}.$ To this end, we have the following conditional result.

\begin{theorem}\label{ttwo}
Let $\tau\gg1$, $\gamma\in [0,1)$, and suppose there exists a convex body $K\in{\cal K}$, such that
\begin{equation}|\{\tau\partial K \cap {\mathbb Z}^d\}| \;\geq\;
C_K\,  \tau^{d-2+\gamma}. \label{points}\end{equation} For any $s\in(0,d)$, there exists a measure $\mu\in{\cal
M}_s$, supported in the unit ball, such that for $p=\frac{d}{s}$ and $t=\tau^\frac{p}{p-1},$ one has
\begin{equation}\label{mest} \sigma_{\mu, K}(t)\;\equiv\;\int_{\partial K} {|\widehat{\mu}(t\omega)|}^2 d\omega_K
\;\geq \;c_K \,t^{-s+(\frac{2s}{d}-1)+ \gamma\frac{p-1}{p}},
\end{equation}
where $d\omega_K$ is the Lebesgue measure on $\partial K$.
\end{theorem}
\begin{corollary}\label{cortwo}
In dimension $2$, there exists $K\in{\cal K},$ such that for a sequence of values of $t$ going to infinity,
there exists a measure $\mu_t\in{\cal M}_s$, supported in the unit ball, such that
\begin{equation}\label{mesttwo} \int_{\partial K} {|\widehat{\mu_t}(t\omega)|}^2 d\omega_K
\;\geq c_K\; t^{-\frac{1}{2}-\frac{s}{4}}.
\end{equation}
\end{corollary}
We remark that (\ref{mest}), even for $\gamma=0$, is always non-trivial for $s>\frac{d}{2}$, while
(\ref{mesttwo}) is non-trivial for $s>\frac{2}{3}$.

\section{Proofs of theorems}
We start out with a simple calculation showing that $\mu_s$ defined
by (\ref{measure}) is in ${\cal M}_s$.
\begin{lemma}
For $s=\frac{d}{p}$, we have $I_s(\mu_s)\lesssim 1$. \label{clc}
\end{lemma}
Clearly, the approximate expression (\ref{mu}) is good enough to substitute for $\mu_s$ in the energy
computation, see (\ref{energy}). For any $x$ in the support of $\mu_s$, let us split
$$\int \frac{d\mu_s(y)}{\|x-y\|^s} = I_1 + I_2,$$ where $I_1$
is taken over the ball $B(x, \frac{c}{q})$ and $I_2$ over its complement. Then
$$ I_1 \lesssim q^{(p-1)d} \int_{B(0,q^{-p})} \frac{dy}{\|y\|^{\frac{d}{p}}} \lesssim 1.$$ Besides, as $A$
is well-distributed, and the $\mu_s$-mass of each peak centered at $a\in q^{-1}A$ in (\ref{mu}) is approximately $q^{-d}$, one has
$$ I_2 \lesssim q^{-d} \sum_{a\in q^{-1}A_q\setminus B(0, cq^{-1}) } \frac{1}{\|a\|^s} =q^{-d+s} \sum_{a\in A_q\setminus B(0, c)} \frac{1}{\|a\|^s}\lesssim 1.$$

\subsubsection*{Proof of Theorem
\ref{tone}} In this section, let us drop the ${\,}_s$ subscript for $\mu_s$ and $\rho_s$ from (\ref{measure}),
to avoid having too many indices. The proof contains three steps, and we start out with two preliminary
observations. Since
\begin{equation}\label{ft}\widehat\rho(\xi)=C_A \widehat\phi(q^{-\frac{d}{s}}\xi)\sum_{a\in q^{-1}A_q} e^{-2\pi i a\cdot\xi},\end{equation} the ``dimension" $s$ characterizing the thickening $\mu_s$ of the atomic measure $\sum_{a\in q^{-1}A_q}\delta(x-a)$ appears only in the cut-off $\widehat\phi(q^{-\frac{d}{s}}\xi).$ Hence, given $s$, it suffices to consider $t=\Theta(q^{\frac{d}{s}})$ only. This is assumed throughout Step 1 of the proof. Indeed, instead of considering $t\ll q^{\frac{d}{s}}$, one can rather increase $s$ (it is assumed that $t\gg q$). In Step 2 we verify the estimate (\ref{sigmabound}) for $s=\frac{d}{2}$ and  $t\lesssim q^2$. Technically, in the end, we will consider separately the ``endpoint case" given $s\in[\frac{d}{2},d)$ and $t= Nq^{\frac{d}{s}},$ for $N\rightarrow\infty$. In this case, $\widehat\phi(q^{-\frac{d}{s}}\xi)$, with $\|\xi\|=t$ satisfies the standard decay estimate $O((1+N)^{-n})$, for any $n$, and this accounts for the pre-factor $\eta$ in the estimates of Theorem
\ref{tone}. This is carried out in Step 3 of the proof.

The second standard preliminary observation is that the density $\rho$ in (\ref{measure}) can be multiplied by any test function that equals one in the unit ball, reflecting the fact that $\mu$ is compact. This implies that
$\widehat{\mu}$ changes slowly on the length scale $\Theta(1)$. Namely, if $A(t,c)$ denotes the spherical shell of
radius $t$ and width $2c$, we have
\begin{equation}\label{ann}\begin{array}{lll}
 t^{d-1}\sigma_{\mu}(t)&\lesssim&\int_{A(t,c)}|\widehat\rho(\xi)|^2 d\xi\\ \hfill
 \\&=& \sup_{f:\;\|f\|_2=1,\;\mbox{supp}\,f\in A(t,c)} \left(\int
 \widehat\rho(\xi) f(\xi)d\xi\right)^2 \\ \hfill \\
 &\equiv & (M_{\mu}[f])^2 \end{array}
\end{equation}

{\sf Step 1.} Take any such $f$. Let $P_q$ be a maximum $\frac{c_1 q}{t}$ separated set on $S^{d-1}$, for some
sufficiently small $c_1$. For $p\in P_q$, let $f_p$ be the restriction of $f$ on the intersection of $A(t,c)$
with the cone, emanating from the origin and built upon the Voronoi cell of $S^{d-1}$ centered at $p$. (The
latter is defined as the set of all points on $S^{d-1}$ that are closer to $p$ than to any other point of
$P_q$.) Decompose\footnote{Technically, one can always smoothen the Voronoi cells out by tweaking them a bit,
see e.g. \cite{T}, but here this is not necessary. Besides, we hope that the fact that the symbol $p$ has
appeared earlier as $p=\frac{d}{s}$ and is used throughout the proof of Theorem \ref{tone} as the summation
index over the partition of $S^{d-1}$ does not cause ambiguity.}
\begin{equation}\label{fdc}f =\sum_{p} f_p,\end{equation} clearly,
\begin{equation}|P_q| \approx \;\left(\frac{t}{c_1 q}\right)^{d-1}\label{many}. \end{equation} By
choosing a small $c_1$, we can ensure that for $c_2$ as small as necessary (in terms of the bounding constants
$c_{\cal A},\,C_{\cal A}$ characterizing the well-distributed set class ${\cal A}$) the support of $f_p$ is
contained in some $d$-dimensional rectangle (henceforth simply {\em rectangle}) of the size $c_2(q^2/t\times
q\times\ldots\times q )$, which is centered at $tp$ and the first measurement is taken in the direction of $p$.

By orthogonality,
\begin{equation}
1=\|f\|_2^2= \sum_{p} \|f_p\|_2^2. \label{ort}\end{equation} Then
\begin{equation}
(M_\mu[f])^2 \leq \left( \sum_p |M_\mu[f_p]| \right)^2.
\end{equation}
To prove (\ref{upper}), we are going to show that for each  $p\in P_q$,
\begin{equation}
|M_\mu[f_p]| \lesssim\; \|f_p\|_2,\;\mbox{ uniformly in } p, \label{prove}\end{equation} as a coarse
estimate. This will imply by Cauchy-Schwartz and (\ref{ort}) that
\begin{equation}
(M_\mu[f])^2\;\lesssim \|f\|_2^2\cdot\sum_{p\in {P_q}}1\lesssim \left(\frac{t}{q}\right)^{d-1}.
\label{pnum}\end{equation}

Let us first prove (\ref{prove}). Without loss of generality, assume that $p=(1,0,\ldots,0)$, relative to the
coordinates $(\xi_1,\xi_2)$, where $\xi_1$ is one-dimensional and $\xi_2$ is $(d-1)$-dimensional.

By Plancherel we have
\begin{equation}
M_\mu[f_p] = \int \widehat{f}_p(x,y) \rho(x,y)dxdy, \label{pl}\end{equation}  where $x$ is one-dimensional
and $y$ is $(d-1)$-dimensional.

The function $f_p$ is supported in the translate by $tp$ of the rectangle $\bar R_p$, where $\bar R_p$ -- with
the above choice of $p=(1,0)$ -- is a ``vertical" rectangle centered at the origin in the $(\xi_1,\xi_2)$
``plane", of width $c_2 q^2/t$ and height (meaning the $\xi_2$-directions) $c_2 q $. Let us write $f_p(\xi)=
h_p(\xi - tp),$ i.e. $h_p$ is supported in $\bar R_p$. All the rectangles involved are further identified with
their characteristic functions.

By the uncertainty principle, as $h_p=h_p \cdot \bar R_p$, its Fourier transform $\widehat h_p$ is approximately
constant in the translates of the dual to $\bar R_p$  rectangle $R_p$ of size $C_2 (t/q^2\times
q^{-1}\times\ldots\times q^{-1})$, relative to the coordinates $(x,y)$. More precisely, if
$\bar{r}_p(\xi_1,\xi_2)$ is a test function which is one in $\bar R_p$ and vanishes outside, say $2\bar R_p$,
then $\widehat h_p=\widehat h_p
* \widehat{\bar{r}}_p$.

Accordingly, let us decompose
\begin{equation} \widehat{h}_p \;=\;\sum_j\widehat{h}_p R_{p,j}\;\equiv\;\sum_j \widehat{h}_{p,j}: \;\;\;\; \|\widehat{h}_p\|^2_2 \;= \;\sum_{j} \|\widehat{h}_{p,j}\|^2_2.\label{ort2}\end{equation} Above,
$R_{p,j}$ are the translates of $R_p$ that together tile some square, covering the unit ball, where $\mu$ is
supported; $R_{p,j}$ is identified with its characteristic function. The constant $C_2$ can be made as large as
necessary by decreasing $c_1$ above. We shall further use well-distributedness of the set $A$, by claiming that
each $R_{p,j}$ supports $q^d\Theta(|R_p|)$ members of $q^{-1}A$.

By Young's inequality
\begin{equation}
\|\widehat{h}_{p,j}\|_{\infty}\;\lesssim\; \frac{1}{\sqrt{|R_p|}}
\|\widehat{h}_{p,j}\|_2, \label{abs}\end{equation} moreover as
$\widehat h_p=\widehat h_p * \widehat{\bar{r}}_p$, we can write
\begin{equation}
\widehat h_{p,j} \;= \; \left(\frac{1}{\sqrt{|R_p|}} \|\widehat{h}_{p,j}\|_2\right)\, R_{p,j}
\psi_{p,j}.\label{rpr}\end{equation} Above, $\psi_{p,j}$ is a smooth function which is $O(1)$ and can be made to
vanish outside $2 R_{p,j}$; in addition one has uniform bounds \begin{equation} |\partial_x
\psi_{p,j}(x,y)|\;=\;O\left(\frac{q^2}{t}\right),\;\;\; |\partial_y \psi_{p,j}(x,y)|
\;=\;O(q).\label{drv}\end{equation}

Clearly
\begin{equation}
\widehat{f}_p \;= \;e^{-2\pi i tx} \widehat{h}_p. \label{fh}\end{equation} I.e. $\widehat{f}_p$ is the rapid
phase $e^{-2\pi i tx}$ that does not depend on $y$ times $\widehat{h}_p$ (this is specific for the spherical
average, versus non-isotropic $\partial K$-averages) which is approximately constant in each rectangle
$R_{p,j}$, with the sharp bound (\ref{abs}).

By (\ref{pl}) we have then

\begin{equation} \begin{array}{llccll}
M_\mu[f_p] & = & \sum_j \int \widehat h_{p} {R_{p,j}}(x,y) e^{-2\pi i tx} \rho(x,y) dy dx \\ \hfill \\ &\equiv &
\sum_j \left( \frac{1}{\sqrt{|R_p|}} \right) \|\widehat{h}_{p,j}\|_2\;\widehat{\tilde\mu}_{p,j}(t),\end{array}
\label{jdecomp}
\end{equation}
where
\begin{equation}\label{muj}\tilde\mu_{p,j}(x) = \int
{R_{p,j}}(x,y)\psi_{p,j}(x,y)\rho(x,y)dy.\end{equation} Now the desired inequality
\begin{equation}\label{pprove} |M_\mu[f_p]|^2\;\lesssim \; \|f_{p}\|_2^2\end{equation}
follows by Cauchy-Schwartz from the trivial bound
\begin{equation} \forall j,\;\;\int R_{p,j}(x,y)\rho(x,y)dxdy \;\lesssim\; |R_p|, \label{coarb}\end{equation} by
well-distributedness of $A$ ($C_2$ has been chosen large enough) and the fact that there are $O(|R_p|^{-1})$
terms in the summation in $j$. This proves (\ref{prove}).

\medskip
{\sf Step 2.} Naturally, cf. (\ref{jdecomp}), similar to
(\ref{muj}), one is tempted to define
\begin{equation}
\mu_{p,j} (x) \;=\; \int R_{p,j}(x,y)\rho(x,y) dy,\label{localmu}\end{equation} and have
$\widehat{\mu}_{p,j}(t)$ substitute  $\widehat{\tilde\mu}_{p,j}(t)$ in the second line of (\ref{jdecomp}). The
two can be related point-wise however only if the $x$-measurement $C_2\frac{t}{q^2}$ of the rectangle $R_p$ is
$\Omega(1)$, to ensure $|\partial_x \psi_{p,j}(x,y)|\;=\;O(1)$ rather than the first bound in (\ref{drv}).

It is easy to achieve this by changing the partition (\ref{fdc}), (\ref{many}) and essentially repeating the
argument up to this point. In this part of the proof, we assume $s=\frac{d}{2}$ and $t\lesssim q^2$. Let us use
a slightly different decomposition of the sphere $S^{d-1}$, with $P_{\sqrt{t}}$ denoting a maximum
$\frac{c_1}{\sqrt{t}}$ separated subset of $S^{d-1}$. Similar to (\ref{fdc}) and (\ref{many}), decompose
$$f=\sum_p f_p,\;\;\;p\in P_{\sqrt{t}},\;\;\; |P_{\sqrt{t}}|\;\lesssim \;t^{\frac{d-1}{2}}.$$ Now $f_p$
is supported inside the rectangle $\bar{R}_p$ of the size $c_2(1\times \sqrt{t}\times\ldots\times\sqrt{t}).$
Accordingly, its dual ${R}_p$ has the size $C_2(1\times
\frac{1}{\sqrt{t}}\times\ldots\times\frac{1}{\sqrt{t}}).$

We repeat the argument from (\ref{ort2}) through (\ref{muj}), with the same notations, relative to the new
partition $\{f_p\}$, only now we can write

\begin{equation} \tilde\mu_{p,j}(x) \;= \;\varphi_{p,j}(x){\mu}_{p,j}(x),\label{newmu}
\end{equation}
for some test function $\varphi_{p,j}(x)$ of a single variable, which is supported on $[-C_2,C_2]$, and is
$O(1)$, together with is derivative. Above, the quantities $\tilde\mu_{p,j}(x)$ and $\mu_{p,j}(x)$ have been
defined respectively by (\ref{muj}) and (\ref{localmu}), only relative to the rectangles $R_{p,j}$ of the size
$C_2(1\times \frac{1}{\sqrt{t}}\times\ldots\times \frac{1}{\sqrt{t}})$, hence the desired properties of
$\varphi_{p,j}(x)$ that arise after integration in $y$ in $(\ref{muj})$, in view of the bound $|\partial_x
\psi_{p,j}(x,y)|\;=\;O(1)$.

Thus we have $\widehat{\tilde\mu}_{p,j} =
\widehat\mu_{p,j}*\widehat\varphi_{p,j}$, and this implies the bound
\begin{equation}|\widehat{\tilde\mu}_{p,j}(t)|\;\lesssim
\;\sup_{\tau}|\widehat\mu_{p,j}(t-\tau)|\,\eta(cC_2\tau),\label{cbd}\end{equation}
where $c$ is independent of the governing constants $c_1,C_2$, and
the quantity $\eta$ has been defined in the statement of Theorem
\ref{tone}.

In view of this, we can give a more refined bound than (\ref{pprove}) following (\ref{jdecomp}). Using
(\ref{cbd}) and the fact that now $|R_p|\approx t^{\frac{1-d}{2}}$, we obtain, essentially repeating the
argument in Step 1,  that
\begin{equation}\begin{array}{llllll}
|M_\mu[f]|^2&\lesssim&\Sigma_{d/2}(t)& \equiv & \sum_{p\in P_{\sqrt{q}}} \left(\frac{1}{|R_p|}\sum_{j}
|\widehat{\tilde\mu}_{p,j}(t)|^2\right) \\ \hfill \\ &&&\lesssim& \sum_{p\in P_{\sqrt{q}}}
\left(t^{\frac{d-1}{2}}\sum_{j} \sup_\tau |\widehat{\mu}_{p,j}(t-\tau)|^2 \,\eta(cC_2\tau)\right).\end{array}
\label{sigma}\end{equation}
 A coarse bound (\ref{coarse}) follows in exactly the same way as (\ref{upper}) on Step 1. I.e. for both
 partitions of $S^{d-1}$, we have
\begin{equation}
|M_\mu[f]|^2\;\lesssim\; \sum_{p} \left(\frac{1}{|R_p|}\sum_{j}
|\widehat{\tilde\mu}_{p,j}(t)|^2\right).\label{gcb}\end{equation}

\medskip
{\sf Step 3.} So far, the bounds (\ref{upper}) -- (\ref{coarse}) of Theorem \ref{tone} have been justified only
for $t\lesssim q^{\frac{d}{s}}$, with $s\in[\frac{d}{2},d)$ on Step 1 and $s=\frac{d}{2}$ on Step 2. Suppose now
that $t\approx Nq^{\frac{d}{s}},$ where $N$ increases. The impact of this shall be compensated by the choice of
the constant $c_1$, increasing the number of Voronoi cells on $S^{d-1}$, to ensure that $C_2$ remains
sufficiently large. Hence, the constants hidden in (\ref{pprove}) as well as in (\ref{sigma}) will increase as
$N^{d-1}$. On the other hand, built into (\ref{ft}), we have the decay of $\widehat\phi(q^{-\frac{d}{s}}\xi)$.
This clearly enables one to use $|R_p|\widehat\phi(q^{-\frac{d}{s}}\xi),$ with $\|\xi\|=t$ as a coarse bound for
$|\widehat{\tilde\mu}_{p,j}(t)|,$ i.e. multiply $|R_p|$ by $C_n (1+N)^{-n}$ for any $n$, and $t\approx
Nq^{\frac{d}{s}}$. This accounts for the presence of the quantity $\eta$ in (\ref{upper}) and (\ref{coarse}) and
completes the proof of Theorem \ref{tone}.

\subsubsection*{Proof of Theorem \ref{ttwo}}
Let us modify the measure $\mu_s$ in (\ref{measure}) slightly, keeping the same notation, with now again
$p=\frac{d}{s}$:
\begin{equation}\label{measuretwo}
d\mu_s(x)\;=\;\rho_s(x)dx, \;\mbox{ with }\;\rho_s(x)\;=\; C_A \,\phi(x)\, q^{-d} \sum_{a\in q^{-1}A}
\phi(a)\phi_{q^p}(x-a).\end{equation} Lemma \ref{clc} clearly remains true, although in comparison with the
expression (\ref{mu}), the pre-factor $\phi(x)$ has enabled to extend the summation over the whole $q^{-1}A$;
besides each peak at $a\in q^{-1}A_q$ has been weighted by $\phi(a)$.

The analog of (\ref{ft}) is now
\begin{equation}\label{fttwo}\widehat\rho_s(\xi)\; = \;C_A\,
\widehat\phi(\xi)* \left(\widehat\phi(q^{-p}\xi)\sum_{a\in q^{-1}A} \phi(a) e^{-2\pi i
a\cdot\xi}\right).\end{equation}

We now consider the special case $A=\Z^d$ and apply the Poisson summation formula to the sum in $a$, which
results in the summation over the dual to $q^{-1}\Z^d$ lattice $q\Z^d,$ at each of whose elements $b$ there sits
a bump $\widehat\phi(\xi-b).$ I.e.
\begin{equation}\label{fttwo1}\widehat\rho_s(\xi)\; \approx \;
\widehat\phi(\xi)* \left[\widehat\phi(q^{-p}\xi)\sum_{b\in q\Z^d} \widehat\phi(\xi-b)\right].\end{equation}

Consider now the average
$$
\int_{\partial K} |\widehat\mu_s(t\omega )|^2 d\omega_K\;\sim\; t^{1-d} \int_{A_K(t,c)} |\widehat\mu_s(\xi )|^2
d\xi,
$$
where $A_K(t,c) = (t+c)K\setminus (t-c)K$. Strictly speaking in the above relation one should have the
$\lesssim$ symbol, however the right-hand side will suffice for the lower bound in this particular case as well.

Indeed, if $K$ is such as stated by Theorem \ref{ttwo}, $t\partial K$ contains
$\Omega\left((\frac{t}{q})^\gamma\right)$ points of the lattice $q\Z^d$, and hence by (\ref{fttwo1}), there are
$\Omega\left((\frac{t}{q})^\gamma\right)$ bumps, each of the hight approximately one and with an
$\Omega(1)$-overlap with the shell $A(t,c)$ or the dilated boundary $t\partial K$ itself.

Therefore,
$$
\int_{\partial K} |\widehat\mu_s(t\omega )|^2 d\omega_K\;\gtrsim\;
t^{1-d}\left(\frac{t}{q}\right)^\gamma,
$$
and the proof of Theorem \ref{ttwo} is complete by choosing
$q=t^{\frac{s}{d}}$, with $\tau =\frac{t}{q}$ in the condition
(\ref{points}).

\medskip
To prove Corollary \ref{cortwo} in the case $d=2$, it is easy to see that some $\tau$-dilates of a piece of the
parabola $\{y=\pm\sqrt{x},\,x\in [0,1]\}$ would contain $\Theta(\tau^{\frac{1}{2}})$ integer points. Indeed, the
dilate in question can be written as $\{(x,\pm\sqrt{x}\sqrt{\tau}),\,x\in[0,\tau]\}$, and if $\tau$ is a square,
the dilate obviously contains an integer point whenever $x$ is a square. The above parabola can be made part of
the boundary $\partial K$ of the body $K$ determining the metric $\|\cdot\|_K$. This proves the corollary.

\begin{remark} Observe that the condition (\ref{points}) can be relaxed by having the points of $q\Z^d$ located
$c$-close to $tK$, rather than immediately on it. In other words, it suffices
to take the right-hand side of (\ref{points}) as the lower bound for
the number of integer points located $c\tau^{\frac{1}{1-p}}$ close
to $\tau
\partial K.$ In the case $p=2$, $s=\frac{d}{2}$, we have
$\tau=q$.\end{remark}

\section{Implications for distance conjectures and lattice point distributions}

Let us first follow up on the discussion surrounding (\ref{fm})--(\ref{ub}) relating the spherical average and
the distance conjectures. Mattila (\cite{Ma}) reformulated the Falconer conjecture as a claim that a compact set
$E$ with $\alpha={\rm dim}_H\,E>\frac{d}{2}$
 should support a Borel probability measure
$\mu\in {\cal M}_s$, for $\frac{d}{2}<s<\alpha,$ such that the corresponding distance measure $\nu_\mu$ has an
$L^2$ density. He then showed that after scaling $\nu_\mu\rightarrow \nu_\mu t^{\frac{1-d}{2}}$, the Hankel
transform
\begin{equation} \widehat{\nu}_\mu(t) \;\equiv\; \int_0^\infty
\sqrt{t\tau}J_{\frac{d}{2}-1}(t\tau)d\nu_\mu(\tau) \;\approx \;t^{\frac{d-1}{2}}\sigma_\mu(t),
\label{ht}\end{equation} where $\sigma_\mu(t)$ is the spherical average (\ref{savr}). Above $J_{\frac{d}{2}-1}$
is the Bessel function of order ${\frac{d}{2}-1}$, and (\ref{ht}) is a variant of the Fourier transform on
$\R_+$, for which the usual properties, such as the the particular Parceval identity, continue to hold. We shall
use the notation $\widehat{\nu}_\mu(t)$ for this (Hankel) transform in the sequel.

Therefore, the Falconer conjecture would follow by the Cauchy-Schwartz inequality if one could bound the second
moment of the distance measure $\nu_\mu$ as
\begin{equation}
\infty\;>\;F(\mu)\;=\;\int_0^{\infty} \sigma^2_\mu(t) t^{d-1} dt
\;\approx \;\|\nu_\mu\|^2_{L^2(\R_+)}. \label{mat}\end{equation}

Note that for $\mu\in {\cal M}_s$, one has the natural energy estimate
\begin{equation}\begin{array}{lllllllll}
\infty & > & I_s(\mu) &=& \int\int
\frac{d\mu_xd\mu_y}{\|x-y\|^s} \\ \hfill \\
&&&\approx & \int \|\xi\|^{s-d}|\widehat{\mu}(\xi)|^2 &\approx&
\int_0^\infty \sigma_\mu(t) t^{s-1}dt,
\end{array}\label{enr}\end{equation} using Plancherel's theorem and then passing to polar
coordinates. Hence, $\sigma_\mu(t)$ is on average $O(t^{-s})$, but
this is not enough for the integral (\ref{mat}) to converge.

The above formalism naturally prompts one to investigate the bounds (\ref{ub}) and define
\begin{equation} \bar\beta(s)\;=\;\sup\{\beta:\,\sigma_\mu(t)\,\leq\, C_{\mu,\beta}
t^{-\beta},\;\;\forall\,\mu \in{\cal M}_s \}.\label{beta}\end{equation} Let also ${\displaystyle
\bar\beta=\limsup_{s\rightarrow\alpha}\bar\beta(s)}$. It follows from (\ref{fm}), (\ref{enr}) that
$\bar\beta\leq\alpha={\rm dim}_H(E).$

Assuming the estimate (\ref{beta}), for any $\mu\in{\cal M}_s$, it follows from (\ref{enr}) that
\begin{equation}F(\mu) \;\leq \; C_\mu \int_0^\infty \sigma_\mu(t) t^{d-\beta-1}dt\;\approx\;
I_{d-\bar\beta(s)}(\mu)\;<\infty,\;\mbox{ for }\;d-\bar\beta(s)\,\leq\, s. \label{lmt}\end{equation} Hence
Falconer conjecture holds if
\begin{equation}\alpha\;>\;d-\bar\beta(\alpha).\label{limit}\end{equation}
An estimate $\bar\beta(s)\geq s-1$ (implicit in \cite{Fa}, \cite{Ma}, and explicit in \cite{Sj}) implies that
the Falconer conjecture is true for $\alpha>\frac{d+1}{2}$, and hence (\ref{beta}) is of major interest for
$s\in[\frac{d}{2},\frac{d+1}{2}]$, as was pointed out earlier.

This formalism extends to the case of $K$-distances, concerning the surface average $\sigma_{\mu,K}(t)$, defined
in (\ref{mest}). Then,  see \cite{HI}, the Mattila formulation of the Falconer conjecture for non-isotropic
distances $\|\cdot\|_K$ is equivalent to proving that
\begin{equation}\label{matk}
F_K(\mu)\;=\;\int_0^\infty \sigma^2_{\mu,K^*}(t)t^{d-1}dt\;<\;\infty,
\end{equation} where
$$
K^*=\{x:\;\sup_{y\in\partial K}x\cdot y\leq1\} $$ is the dual body of $K\in {\cal K}$.

\subsubsection*{Proof of Corollary \ref{corone}}
Theorem \ref{tone} implies that the measure $\mu_s$ defined by (\ref{measure}) satisfies (\ref{ub}) for
$\beta=\frac{d-1}{d}s$. Hence, by (\ref{lmt}), the Falconer conjecture is satisfied by the support of
$\mu_s,$ provided that $s\geq \frac{d^2}{2d-1}$. Therefore, the number of distinct $q^{-p+1}$ separated distances generated by the set $A_q$, where $p=\frac{d}{s}$, is bounded from below by a constant times
\begin{equation} q^p = q^{\frac{2d-1}{d}}\approx |A_q|^{\frac{2}{d}-\frac{1}{d^2}}. \label{dst}\end{equation}

This proves Corollary \ref{corone}. Let us point out here that this is precisely the lower bound obtained by
Moser (\cite{Mo}) in the case $d=2$ (see also \cite{Io} for higher-dimensional generalization of this method),
and Solymosi and Vu (\cite{SV}) for well-distributed sets using methods of geometric combinatorics. Recently
Solymosi and T\'oth (\cite{ST}) made further progress in that direction, having improved the margin
$\frac{1}{d^2}$ in (\ref{dst}) to $\frac{2}{d(d^2+1)}$.

\subsubsection*{General bounds for the spherical average}
The strongest general spherical average bounds summarized in (\ref{gnr}) for $s\in[\frac{d}{2},\frac{d+1}{2}]$
are due to Wolff (\cite{WI}) in the case $d=2$ and Erdo\~gan (\cite{Er}) in higher dimensions. More precisely,
Wolff showed that a general $\mu\in {\cal M}_s$, for $s\geq 1$, satisfies (\ref{ub}) with any
$\beta<\frac{s}{2}.$ This cannot be improved beyond the endpoint, because Sj\"olin (\cite{Sj}) used a Knapp-type
example to show that for $s\geq 1,$ there are measures in ${\cal M}_s$ that satisfy
\begin{equation}
\sigma_\mu(t) \gtrsim t^{-\left(\frac{s}{2}+\frac{d-2}{2}\right)} I_s(\mu). \label{sjolin}\end{equation}
Namely, one has
\begin{equation}
\label{wolff}\bar\beta(s)=\frac{s}{2}, \mbox{ for } s\geq 1 \mbox{ and } d=2.
\end{equation}

Observe on the other hand that as $\bar\beta(s) \leq \alpha$, the estimate (\ref{sjolin}) provides non-trivial information only in the range of Hausdorff dimensions $\alpha>d-2$.

Sj\"olin's example shows that in dimensions $2$ and $3$, the Falconer conjecture cannot be resolved in full generality merely by proving sharp power majorants for the spherical average $\sigma_\mu(t)$, but leaves open the question whether this may be possible for $d\geq4$. This question has been recently asked by Erdo\~{g}an (\cite{Erg}) who generalized Wolff's result to $d\geq 2,$ obtaining the best known upper bounds for
$\bar\beta(s)$ in higher dimensions, although not necessarily unimprovable in higher dimensions, cf.
(\ref{gnr}):
\begin{equation} \bar\beta(s) \geq \frac{d+2s-2}{4},\mbox{ for } s\in\left[\frac{d}{2}, 1+\frac{d}{2}\right].\label{cbb}\end{equation} Our bound (\ref{cb}) is an improvement
over (\ref{cbb}), and it appears reasonable to ask the following.
\begin{question}\label{q1} Does (\ref{cb}) generalize to the class ${\cal M}_s$ (at least in the important case
$s\in[\frac{d}{2},\frac{d+1}{2}]$) in the case $d\geq 3$, and if it does, is it generally best
possible?\end{question} We believe that the first part of the question can be answered affirmatively. Observe
that the proof of Theorem \ref{tone} is still valid if we consider a general $\mu\in{\cal M}_s$, which is
well-distributed in the sense that for some $q\gg1$, any ball $B(x,q^{-1})$ of radius $\Omega(q)$ has the
$\mu$-mass of approximately $q^{-d}$, i.e. on the length scales $q^{-1}$ and above, $\mu$ approximates the
Lebesgue measure. Let us give some heuristics when this is not the case, i.e., the parameter $q$ is not built
into the problem. Then, in order to evaluate the spherical average $\sigma_\mu(t)$, one can effectively
(eventually losing the endpoint, due to the issues of dimension) assume that $\mu$ is a density supported on a
union of disjoint balls of radius $t^{-1},$ so that the $\mu$-mass of each ball is $O(t^{-s})$. As $\mu$ is a
compactly supported probability measure, the total number of such balls in its support is approximately $t^s$,
and $q=t^{\frac{s}{d}}$ arises as a natural partition parameter, in the sense that a subcube of diameter
$q^{-1}$ contains on average one of the union of balls whereupon $\mu$ is supported. Therefore, the partition of
the sphere of radius $t$ onto pieces of diameter $O(q)$ arises naturally, and would yield the analog of the
double sum given by (\ref{gcb}). The trivial estimate that has been applied to the double sum, which claimed
that the expression in brackets there was $O(1)$ for each partition angle $p$, is no longer applicable, because
given $p$, the mass $\mu$ may not be distributed between the tiles $R_{p,j}$ uniformly. Then one has to tackle
the whole double sum in (\ref{gcb}). This would create a reasonably accessible combinatorial problem, similar to
the one underlying the general proof of the two-dimensional case in \cite{WI}, as well as proofs of the recent
sharp bilinear restriction theorems of Wolff (\cite{WA}) and Tao (\cite{T}). The results of the latter paper
were adapted to estimate the spherical average by Erdo\~{g}an (\cite{Erg}), with no evidence of being sharp in
$d\geq3$. It is important that in any case such an approach would still ignore the phases in (\ref{jdecomp}), so
the estimates (\ref{cb}) as well as (\ref{cbb}) are in essence {\em coarse estimates}.

We see no cogent reason to believe or disbelieve whether the bound (\ref{cb}) may be tight in higher dimensions,
for the lack of geometric concept that would underly a possible counterexample. Note that Theorem \ref{ttwo}
falls short of doing this, for in the spherical, alias Euclidean case one has $\gamma=0$. The estimate
(\ref{mest}) implies that
\begin{equation}
\bar\beta(s)\;<\;s,\;\;\mbox{ for }\;\;s>\frac{d}{2}, \label{hd}\end{equation} leaving the critical case
$s=\frac{d}{2}$ open. On the other hand, from the point of view of this paper, Sj\"olin's example can be
rendered essentially one-dimensional. Namely, the discretized version thereof is as follows. In the plane, one
takes points with coordinates $(x,y)=(j/q,0),\,j=0,\pm1,\ldots,\pm q$, thickens them into rectangles of width
$q^{-\frac{d}{s}}$ and height $q^{-1},$ puts a uniform probability measure thereon and uses the one-dimensional
Poisson summation formula to look at the Fourier side at $t\approx q^{\frac{d}{s}}$. I.e., within the
decomposition framework of Theorem \ref{tone}, it is basically equivalent to just having a single direction $p$
in the double sum (\ref{gcb}). More precisely, other values of $p\in S^1$ do not contribute due to cancelations.
This does not suffice to match the upper bound (\ref{cb}) in the case $d=3$, and is even less potent in higher
dimensions. It also indicates that cancelations for different values of $p$ are inherent in the problem. More
precisely, the measures $\mu_{p,j}(x)$, defined by (\ref{localmu}), and localizing $\mu$ in different directions
$p$, cannot all resonate with fast plane waves in these directions, all having the same frequency $t$.  Hence
targeting the sharp bounds for the sum in (\ref{gcb}), one cannot simply ignore the presence of the phase
factors in (\ref{jdecomp}).

The estimate (\ref{sigmabound}) of Theorem \ref{tone} is conditional on the term  $\Sigma_{d/2}(t)$ which is
given explicitly by (\ref{sigma}). Naturally, the estimate  (\ref{sigmabound}) poses a  question of estimating
the quantity $\Sigma_{d/2}(t)$ by using something more intricate than coarse estimates. Observe that if
$s=\frac{d}{2}$, then $t\approx q^2$, and the underlying well-distributed set $A$ is the lattice $\Z^d,$ the
number of nonzero terms in the summation over $p$ in (\ref{sigma}) will be approximately $q^{d-2}$ (modulo a
slowly growing function of $q$ in the case $d=2$). Indeed, any such $p$ would correspond to a point of the
lattice $q\Z^d$ lying in the $O(1)$ neighborhood of the sphere of radius $t\approx q^2.$ The number of such
points cannot exceed $O(q^{d-2})$ in dimensions three and higher, with an additional slowly growing term in
dimension two. This implies that in this specific case, for $s=\frac{d}{2}$, the bound (\ref{sigma}) improves
from $t^{-\frac{d-1}{2}}$ to $t^{-\frac{d}{2}}$ (modulo a slowly growing function of $q$ in the case $d=2$)
which is precisely what one needs to prove the Erd\"os/Falconer conjecture. As we mentioned earlier, it was the
lattice example that inspired the distance conjectures (\ref{ec}) and (\ref{fc}). The following question
continues in this vein.
\begin{question}\label{q2} Is it true that for the measures (\ref{measure}) on thickenings of well-distributed sets and $s=\frac{d}{2}$, one actually has
\begin{equation} \sigma_{\mu_{d/2}(t)} \lessapprox t^{-\frac{d}{2}}? \label{cj}\end{equation}\end{question} This cannot be true for $s>\frac{d}{2}$ by (\ref{hd}).
Besides, the answer may possibly be positive only for the Euclidean, spherical average, since generalizing the
spherical average to $\sigma_{\mu_s,K}$ is impossible in light of Corollary (\ref{cortwo}). By (\ref{limit}),
the affirmative answer would also imply the Erd\"os conjecture for well-distributed sets. In this sense, as a
special feature of the sphere, this question is analogous with the Erd\"os single distance conjecture mentioned
further below.

\subsubsection*{Non-isotropic surface means and lattice point distributions} Theorem \ref{ttwo} implies
that the Falconer conjecture for a class of distances $\|\cdot\|_{K\in{\cal K}}$ cannot be resolved by proving
best possible majorants for the quantity $\sigma_{\mu,K}$ alone, provided that for any fixed $\varepsilon>0$ and
arbitrarily large $\tau,$ there exists $K\in{\cal K}$ such that the condition (\ref{points}) holds for any
$\gamma>0.$ Let us address the issue how large $\gamma$ can possibly be. Much more is known to this effect in
$d=2$ versus higher dimensions.

In the case $d=2$, by the result of Bombieri and Pila (\cite{BP})  there are no $C^\infty$ bodies $K$, such that
(\ref{points}) holds for $\gamma=\frac{1}{2}+\epsilon,$ for any $\epsilon>0$. (More recently, the result of
Bombieri and Pila has been given more refinement under additional assumptions which are beyond the scope of this
paper.) The conjecture of Schmidt (\cite{Sch}) states that this is actually the case for the class $C^2$.
Observe that so far, for the analysis of $\partial K$ means in the literature dedicated to the Falconer
conjecture a finite order of differentiability of $\partial K$ suffices.

In higher dimensions, to our best knowledge, there are no explicit examples of  $K$ satisfying (\ref{points})
with $\gamma>0$. An upper bound for $\gamma$ can be derived, for instance, from the results concerning the
lattice point distributions error term. We now quote the estimate due to M\"uller (\cite{Mu}). Let
$N(\tau)=|\{\tau K \cap {\mathbb Z}^d\}|$, $d > 2$, suppose $N(\tau)=|K|\tau^d+E(\tau)$. Then
\begin{equation} |E(\tau)| \lesssim \tau^{d-2+\gamma_d},\mbox{ with } \gamma_d =\left\{ \begin{array}{lll}\frac{20}{43}, & d=3, \\ \hfill \\ \frac{d+4}{d^2+d+2}, & d \ge 4. \end{array} \right.\label{muller}\end{equation} Clearly (\ref{muller}) implies $\gamma\leq\gamma_d$ for the condition (\ref{points})
for otherwise one could construct an immediate counterexample to (\ref{muller}).

Returning to $L^2$ surface averages, the above quoted upper bounds (\ref{wolff}) and (\ref{cbb}) of Wolff and
Erdo\"gan are applicable to the quantity $\sigma_{\mu,K}$ defined in (\ref{mest}) as well, with the bounding
constants now depending on $K$. In the same fashion, one can easily see from proof of Theorem \ref{tone} that
the coarse bounds (\ref{upper}) and (\ref{cb}) are also applicable in this case. If one attempts to use these
bounds for the specific case $A=\Z^d$, they lead to a trivial estimate $\gamma\leq1$. In other words,
$\tau\partial K$ contains no more than $\tau$ integer points thereon or in its $\frac{1}{\tau}$-vicinity.

Observe that the affirmative answer to Question \ref{q2} would imply that $\tau S^{d-1}$ contains no more than
$C_\epsilon \tau^{\epsilon}$ integer points, for any $\epsilon>0$, which is indeed known to be true. Such an
improvement could in principle come from taking the phase factors in (\ref{jdecomp}), (\ref{sigma}) into
account. The fact that these phase factors appear in their present form is the special feature of the Euclidean
case. We have the following generalization of Question \ref{q2}.
\begin{question}\label{q3} Are the bounds for the number of lattice points on or near the dilates of the
boundaries of $K\in{\cal K}$ a particular case of general asymptotic bounds for the quantity
$\sigma_{\mu,K}(t)$, for measures arising as thickenings of well-distributed sets and not necessarily lattices?
Is it true in particular that in the case $d=2$, cf. (\ref{mu}),
\begin{equation}\sigma_{\mu_{d/2},K}(t)\;\lessapprox \; t^{-\frac{3}{4}}?\label{tq}\end{equation}\end{question}

\subsubsection*{Single distance conjecture and the spherical average}
We finish with some remarks about the Erd\"os single distance conjecture in the case $d=2$, which in the
well-distributed setting can be written as \begin{equation}\label{wdse}\sup_{u\in\Delta(A_q)}|\{(x,y)\in
A_q\times A_q: \|x-y\|=u\}|\;\lessapprox q^2.\end{equation} Using (\ref{st}), with $d=2$ and $s=\frac{d}{2}=1$,
it follows that in terms of the measure (\ref{measure}) and its spherical average (\ref{savr}) it is equivalent
to asking for any $\tau\in(0,1)$ that
\begin{equation}\label{wdseone}\int_0^{\infty} t J_0(\tau t)\sigma_{\mu_{d/2}}(t)dt
\lessapprox 1.\end{equation} At the same time, the Mattila criterion (\ref{mat}) is
\begin{equation}\label{wdsetwo}\int_0^{\infty} t \sigma^2_{\mu_{d/2}}(t)dt \lessapprox 1.
\end{equation} Note that by (\ref{ft}), it is essentially sufficient to integrate up to $q^2$.
The affirmative answer to Question \ref{q2}, the bound $\sigma_{\mu_{d/2}}(t)\,\lessapprox\,\frac{1}{t}, d=2,$ implies (\ref{wdsetwo}), but not (\ref{wdseone}) which requires more regularity than merely the above majorant for the spherical average $\sigma_{\mu_{d/2}}(t)$. By using the asymptotics of the Bessel function $J_0$, the well-distributed set single distance conjecture reduces to the estimate
\begin{equation}\int_{q^2}^{2q^2} \sqrt{t} e^{-2\pi i \tau t}
\sigma_{\mu_{d/2}}(t)dt\;\lessapprox 1,\forall \tau \in(0,1). \label{more}\end{equation} Hence, (\ref{more})
asserts a special property of the Euclidean distance, which is more stringent than (\ref{cj}) in Question
\ref{q2} (in the case $d=2$). This certainly adds to the credibility of the conjectured affirmative answer to
the latter question.

Observe that the single distance conjecture, as well as (\ref{cj}),  are generally not true in the distance
class $\|\cdot\|_{K\in{\cal K}},$ the counterexample to the former conjecture being constructed in essentially
the same way as it has been done to prove Corollary \ref{cortwo}.

The best known single distance conjecture bound of the form (\ref{wdse}) is $q^{\frac{8}{3}}$ and is due to
Spencer, Szemer\'edi, and Trotter (\cite{SST}). It arises as an immediate corollary of the Szemer\'edi-Trotter
theorem. In light of the discussion in this paper, the latter theorem can be formulated as follows. Given a set
$P$ of points and the set $J$ of translates $T_{j\in J}$ of $K$, with $|P|\approx|J|=q^2$, one has the following
bound for the number of incidences:
\begin{equation}
|\{(p,j)\in P\times J: p\in T_j\partial K\}| \lesssim q^{\frac{8}{3}}. \label{st}\end{equation} The single
distance conjecture then claims that the bound in (\ref{st}) can be improved to $q^2$ (modulo a slowly growing
function of $q$) in the special case $\partial K=S^1$. If one believes in optimality of the parabola example
used to prove Corollary \ref{cortwo}, the final (and perhaps the most difficult) question we ask in this section
is the following.
\begin{question}\label{q4} Is this true that for  a general $K\in{\cal K}$, the bound $q^{\frac{8}{3}}$ in (\ref{st}) can be improved to $q^{\frac{5}{2}}$ (modulo a slowly growing function of $q$)?
\end{question}
The parabola example mentioned above shows that the exponent $\frac{5}{2}$ cannot be improved.

\end{document}